\documentclass[leqno]{article}

\newcounter{supersection}[section]
\newtheorem{pr}[supersection]{Proposition}
\newtheorem{thm}[supersection]{Theorem}
\newtheorem{lm}[supersection]{Lemma}
\newtheorem{re}[supersection]{Remark}
\newtheorem{co}[supersection]{Corollary}
\newtheorem{ex}[supersection]{Example}

\renewcommand{\thefootnote}{}

\def\bibname{\textbf{REFERENCES}}
\def\thebibliography#1{\paragraph*{\uppercase{\bibname}}\list
{[\arabic{enumi}]}{\settowidth\labelwidth{[#1]}\leftmargin\labelwidth
\advance\leftmargin\labelsep\usecounter{enumi}}
\def\newblock{\hskip .11em plus .33em minus .07em}
\sloppy\clubpenalty4000\widowpenalty4000
\sfcode`\.=1000\relax}

\def\Dj{D{\hspace{-.75em}\raisebox{.3ex}{-}\hspace{.4em}}}

\def\stop{\mbox{\footnotesize {\vrule width 6pt height 6pt}}}

\renewcommand{\thefootnote}{\arabic{footnote}}

\begin{document}

\centerline{\large \bf ANALYTICAL AND DIFFERENTIAL-ALGEBRAIC}

\medskip

\centerline{\large \bf PROPERTIES OF GAMMA FUNCTION}

\bigskip

\bigskip

\centerline{\bf \v{Z}arko Mijajlovi\'{c}\mbox{${}^{1)}$},
Branko Male\v sevi\'{c}\mbox{${}^{2)}$}}

\bigskip

\begin{center}
{${}^{1)}$Faculty of Mathematics, University of Belgrade, \\
Studentski trg 16, 11000 Belgrade, Serbia and Montenegro  \\[2.0 ex]

${}^{2)}$Faculty of Electrical Engineering, University of Belgrade, \\
Bulevar Kralja Aleksandra 73, 11000 Belgrade, Serbia and Montenegro}
\end{center}

\smallskip

\begin{abstract}
In this paper we consider some analytical relations between gamma function
$\Gamma(z)$ and related functions such as the Kurepa's function $K(z)$ and
alternating Kurepa's function $A(z)$. It is well-known in the physics that
the Casimir energy is defined by the principal part of the Riemann function
$\zeta(z)$ $($Blau, Visser, Wipf; Elizalde$)$. Analogously, we consider
the principal parts for functions $\Gamma(z)$, $K(z)$, $A(z)$ and we also
define and consider the principal part for arbitrary meromorphic functions.
     Next, in this paper we consider some differential-algebraic $($d.a.$)$
properties of functions $\Gamma(z)$, $\zeta(z)$, $K(z)$, $A(z)$. As it is
well-known $($H\" older; Ostrowski$)$ $\Gamma(z)$ is not a solution of any
d.a. equation. It appears that this property of $\Gamma(z)$ is universal.
Namely, a large class of solutions of functional differential equations
also has that property. Proof of these facts is reduced, by the use of the
theory of differential algebraic fields $($Ritt; Kaplansky; Kolchin$)$,
to the d.a. transcendency~of~$\Gamma(z)$.
\end{abstract}

\setcounter{footnote}{1}
{}\footnotetext{{\it Email address}\,:\, {\tt zarkom@eunet.yu}}
\setcounter{footnote}{2}
{}\footnotetext{{\it Email address}\,:\, {\tt malesh@eunet.yu}}

\renewcommand{\thefootnote}{}

{}\footnotetext{Second author supported in part by the project
MNTRS, Grant No. ON144020.}


\section{Analytical properties}

In this section we consider analytical properties of the gamma and related
functions which pertain to the principal part of a function at a point.

\subsection{The principal part of the gamma function}

Gamma function is defined by the integral:
\begin{equation}
\Gamma(z) = \displaystyle\int\limits_{0}^{\infty}{e^{-t} t^{z-1} \, dt},
\end{equation}
which converges for  $Re(z) > 0$.  It is possible to
form analytical continuation of this function over the whole set
of the complex numbers $\mathbf{C}$ except at $z = -k$, where
\mbox{$k \in \mathbf{N}_{0} = \{0, 1, 2, \ldots\}$}. One
approach to analytical continuation is given by:
\begin{equation}
\Gamma(z)
=
\displaystyle\sum\limits_{n=0}^{\infty}{\displaystyle\frac{(-1)^{n}}{n!(n+z)}}
+
\displaystyle\int\limits_{1}^{\infty}{e^{-t}t^{z-1} \, dt},
\end{equation}
see \cite{Markushevich:65}, \cite{Leibbrandt:75}. Residue at $z = -k$, $k \in \mathbf{N}_{0}$, is:
\begin{equation}
\mathop{\mbox{\rm res}}\limits_{z = -k}{\Gamma(z)} = \displaystyle\frac{(-1)^{k}}{k!} \, .
\end{equation}
It is possible to extend the domain of the gamma function to the
set of all complex numbers $\mathbf{C}$ in the sense of the
principal part at a point as follows. For a meromorphic function
$f(z)$, on the basis of Cauchy's integral formula, we define {\em
the principal part at point $a$} (see: \cite{Slavic:70}, \cite{Malesevic:03}):
\begin{equation}
\label{GAMMA_PV_a}
\mathop{\mbox{\rm p.p.}}\limits_{z = a}{f(z)}
=
\lim\limits_{\rho \rightarrow 0_{+}}{\displaystyle\frac{1}{2 \pi i}
\!\!\!\!\displaystyle\oint\limits_{|z-a|=\rho}{\!\!\!\!\displaystyle\frac{f(z)}{z-a}\,dz}}.
\end{equation}
If the point $a$ is regular for the function $f(z)$ then
$\mathop{\mbox{\rm p.p.}}\limits_{z = a}{f(z)} = f(a)$; otherwise
the principal part at the pole $z=a$ exists as a finite complex
number:
\begin{equation}
\mathop{\mbox{\rm p.p.}}\limits_{z = a}{f(z)}
=
\mathop{\mbox{\rm res}}\limits_{z=a}{{\Big (}\displaystyle\frac{f(z)}{z-a}}{\Big )},
\end{equation}
as cited in \cite{Malesevic:06}. Let us determine basic properties
of the principal part at the point. For two meromorphic functions
$f_1(z)$ and $f_2(z)$ additivity holds~\cite{Slavic:70}:
\begin{equation}
\mathop{\mbox{\rm p.p.}}\limits_{z = a}{{\Big (}f_{1}(z)+f_{2}(z){\Big )}}
=
\mathop{\mbox{\rm p.p.}}\limits_{z = a}{f_{1}(z)}
+
\mathop{\mbox{\rm p.p.}}\limits_{z = a}{f_{2}(z)}.
\end{equation}
In the paper \cite{Slavic:70} it is proved that multiplicativity
of the principal part does not hold. Namely for the principal part
the following statement is true.

\begin{thm}
Let $f_{1}(z)$ be a holomorphic function at the point $a$ and let $f_{2}(z)$ be a meromorphic function with
pole of the $m$-th order at the same point $a$. Then$:$
\begin{equation}
\mathop{\mbox{\rm p.p.}}\limits_{z = a}{ {\Big (} f_{1}(z) \cdot f_{2}(z) {\Big )} }
=
\displaystyle\sum\limits_{k=0}^{m}{\displaystyle\frac{f_{1}^{(k)}(a)}{k!}
\mathop{\mbox{\rm p.p.}}\limits_{z = a}{ {\Big (} (z-a)^k \cdot f_{2}(z) } {\Big )} }.
\end{equation}
\end{thm}

\noindent {\bf Proof} \enskip Let $f_{1}(z)$ and let $f_{2}(z)$ be
represented by the series:
\begin{equation}
f_{1}(z) = \displaystyle\sum\limits_{i= 0}^{\infty}{\frac{f_{1}^{(i)}(a)}{i!}(z-a)^{i}}
\quad \mbox{and} \quad
f_{2}(z)
=
\displaystyle\sum\limits_{j= -m}^{\infty}{c_{j}(z-a)^{j}},
\end{equation}
for some $c_{j} \in \mathbf{C}$ $(j \geq -m)$, $c_{-m} \neq 0$ and $z \neq a$.
Let us notice that $\mathop{\mbox{\rm p.p.}}\limits_{z = a}{f_{2}(z)} = c_{0}$.
Multiplying the following series:
\begin{equation}
\displaystyle\frac{f_{1}(z) - f_{1}(a)}{z - a} \cdot f_{2}(z)
\, = \,
\displaystyle\sum\limits_{i=1}^{\infty}{\displaystyle\frac{f_{1}^{(i)}(a)}{i!}(z-a)^{i-1}}
\cdot \! \displaystyle\sum\limits_{j=-m}^{\infty}{c_{j}(z-a)^{j}}
\end{equation}
we obtain:
\begin{equation}
\mathop{\mbox{\rm res}}\limits_{z = a}{ {\Big (} \displaystyle\frac{f_{1}(z) - f_{1}(a)}{z - a} \cdot f_{2}(z) {\Big )} }
=
\displaystyle\sum\limits_{k=1}^{m}{\displaystyle\frac{f_{1}^{(k)}(a)}{k!} \cdot c_{-k}} \; .
\end{equation}
Hence:
\begin{equation}
\begin{array}{rcl}
\mathop{\mbox{\rm p.p.}}\limits_{z = a}{ {\Big (} f_{1}(z) \cdot f_{2}(z) {\Big )} }
& \!\!=\!\! &
\mathop{\mbox{\rm res}}\limits_{z = a}{ {\Big (} \displaystyle\frac{f_{1}(z) \cdot f_{2}(z)}{z-a} {\Big )}}   \\[3.0 ex]
& \!\!=\!\! &
\mathop{\mbox{\rm res}}\limits_{z = a}{ {\Big (} \displaystyle\frac{f_{2}(z)}{z-a} {\Big )}} \cdot f_{1}(a)
+
\mathop{\mbox{\rm res}}\limits_{z = a}{ {\Big (} \displaystyle\frac{f_{1}(z) - f_{1}(a)}{z - a} \cdot f_{2}(z) {\Big )}}
                                                                                                              \\[3.0 ex]
& \!\!=\!\! &
f_{1}(a) \cdot c_{0}
+
\displaystyle\sum\limits_{k=1}^{m}{\displaystyle\frac{f_{1}^{(k)}(a)}{k!} \cdot c_{-k}}                       \\[3.00 ex]
& \!\!=\!\! &
\displaystyle\sum\limits_{k=0}^{m}{\displaystyle\frac{f_{1}^{(k)}(a)}{k!}
\mathop{\mbox{\rm p.p.}}\limits_{z = a}{ {\Big (} (z-a)^k \cdot f_{2}(z) } {\Big )} }. \;\; \stop
\end{array}
\end{equation}
\begin{re}
The phrase "function is holomorphic at the point $a$" means not just function is differentiable at $a$,
but differentiable everywhere within some open disk centered at $a$ in the complex plane.
\end{re}
\begin{co}
Let $f_{1}(z)$ be a holomorphic function at the point $a$ and let $f_{2}(z)$ be a meromorphic function with
simple pole at the same point $a$. Then$:$
\begin{equation}
\mathop{\mbox{\rm p.p.}}\limits_{z = a}{{\Big (}f_{1}(z) \cdot f_{2}(z){\Big )}}
=
f_{1}(a) \cdot \mathop{\mbox{\rm p.p.}}\limits_{z = a}{f_{2}(z)}
+
{f_{1}\!}^{'}\!(a) \cdot \mathop{\mbox{\rm res}}\limits_{z = a}{f_{2}(z)}.
\end{equation}
The previous formula, in the case of the zeta function $f_2(z) = \zeta(z)$,
is also given~in {\rm \cite{BlauVisserWipf:88a}}, {\rm \cite{BlauVisserWipf:88b}}.
\end{co}
For meromorphic function $f(z)$ with simple pole at the point $z = a$ the following formula is true \cite{Malesevic:03}:
\begin{equation}
\label{PV_f_1}
\mathop{\mbox{\rm p.p.}}\limits_{z = a}{f(z)}
=
\lim\limits_{\varepsilon \rightarrow 0}{
\displaystyle\frac{f(a-\varepsilon)+f(a+\varepsilon)}{2}}.
\end{equation}
Especially for gamma function $\Gamma(z)$ it is true \cite{Slavic:70}, \cite{Malesevic:03}:

\vspace*{-2.5 mm}

\begin{equation}
\label{PV_1}
\mathop{\mbox{\rm p.p.}}\limits_{z = -n}{\!\!\Gamma(z)}
=
(-1)^{n}\displaystyle\frac{\Gamma^{'}(n+1)}{\Gamma(n+1)^2}
=
\displaystyle\frac{\mbox{\small $-\gamma$} + \mbox{\footnotesize $\displaystyle\sum_{k=1}^{n}{\displaystyle\frac{1}{k}}$}}{n!},
\end{equation}

\smallskip\noindent
where $\gamma$  is Euler's constant and $n \in \mathbf{N}_{0}$.

\subsection{The principal part of the Kurepa's functions}

{\Dj }. Kurepa introduced in paper \cite{Kurepa:71} function $K(z)$ by integral:
\begin{equation}
\label{K_INT_1}
K(z)
=
\displaystyle\int\limits_{0}^{\infty}{
e^{-t} \displaystyle\frac{t^{z}-1}{t-1} \: dt},
\end{equation}
which converges for $Re(z) > 0$, and it represents one
analytical extension of the sum of factorials:
\begin{equation}
\label{K_SUM_1}
K(n) = \displaystyle\sum\limits_{i=0}^{n-1}{i!} \, .
\end{equation}
For the function $K(z)$ we use the term {\em Kurepa's function} and it is one solution of the
functional equation:
\begin{equation}
\label{K_FE_1}
K(z) - K(z-1) = \Gamma(z).
\end{equation}
Let us observe that it is possible to make analytical continuation
of Kurepa's function $K(z)$ for $Re(z) \leq 0$. In
that way, the Kurepa's function $K(z)$ is a meromorphic function
with simple poles at $z=-1$ and $z=-n$ $(n \!\geq\! 3)$. At point
\mbox{$z = -2$} Kurepa's function has a removable singularity and
$K(-2) \stackrel{\mbox{\tiny def}}{=} \lim\limits_{z \rightarrow
-2}{K(z)} = 1$. Kurepa's function has the following residues:
\begin{equation}
\mathop{\mbox{\rm res}}\limits_{z = -1}{K(z)} = -1
\quad\mbox{and}\quad
\mathop{\mbox{\rm res}}\limits_{z = -n}{K(z)} =
\displaystyle\sum\limits_{k=2}^{n-1}{\displaystyle\frac{(-1)^{k-1}}{k!}}
\quad (n\!\geq\!3).
\end{equation}
Previous results for Kurepa's function are given according to \cite{Kurepa:73} and \cite{Slavic:73}.
The functional equation (\ref{K_FE_1}), besides Kurepa's function $K(z)$, has another solution by series:
\begin{equation}
\label{Def_K1}
K_{1}(z) = \displaystyle\sum\limits_{n=0}^{\infty}{\Gamma(z-n)},
\end{equation}
which converges over the set $\mathbf{C} \backslash \mathbf{Z}$ \cite{Malesevic:03}.

\medskip
\noindent
Extension of domain of functions $K(z)$ and $K_{1}(z)$ in the sense of the principal
part at the point is given by the following statements \cite{Slavic:73}, \cite{Malesevic:03}.
\begin{lm}
Let us define $ L_{1} = -\displaystyle\sum\limits_{n=0}^{\infty}{
\!\,\mathop{\mbox{\rm p.p.}}\limits_{z = -n}{\!\,\Gamma(z)} } $.
Then$:$
\begin{equation}
\label{Const_L_1}
L_{1}
=
\displaystyle\frac{1}{e} {\Big (} \gamma + \displaystyle\sum\limits_{n=1}^{\infty}{\displaystyle\frac{1}{n!n}} {\Big )}
=
\displaystyle\frac{\mbox{\rm Ei}(1)}{e}
\approx
0.697 \, 174 \, 883 \, ,
\end{equation}
where $\mbox{\rm Ei}$ is function of exponential integral.
\end{lm}
\begin{thm}
For the functions $K(z)$ and $K_{1}(z)$ are true$:$
\begin{equation}
\label{K_PV}
\mathop{\mbox{\rm p.p.}}\limits_{z = -n}{\!\!K(z)}
=
-\displaystyle\sum\limits_{i=0}^{n-1}{
\mathop{\mbox{\rm p.p.}}\limits_{z=-i}{\Gamma(z)}}
=
\displaystyle\sum\limits_{i=0}^{n-1}{(-1)^{i+1}
\displaystyle\frac{\Gamma^{'}(i+1)}{\Gamma(i+1)^2}}
\quad (n \!\in\! \mathbf{N})
\end{equation}
and
\begin{equation}
\mathop{\mbox{\rm p.p.}}\limits_{z = n}{K_{1}(z)}
=
\mathop{\mbox{\rm p.p.}}\limits_{z = n}{K(z)} - L_{1}
\quad (n \!\in\! \mathbf{Z}).
\end{equation}
\end{thm}
The connection between functions $K(z)$ and $K_{1}(z)$ is given by Slavi\'{c}'s formula which
is presented in the following statement \cite{Slavic:73}, \cite{Marichev:83}, \cite{Malesevic:03}.
\begin{thm}
It is true$:$
\begin{equation}
\label{K_sum_Slavic}
K(z)
=
\displaystyle\frac{1}{e} {\Big (} \gamma + \displaystyle\sum\limits_{n=1}^{\infty}{\displaystyle\frac{1}{n!n}} {\Big )}
-
\displaystyle\frac{\pi}{e} \, \mbox{\rm ctg} \pi z
+
\displaystyle\sum\limits_{n=0}^{\infty}{\Gamma(z-n)},
\end{equation}
where the values in the previous formula, in integer points $z$,
are determined in the sense of the principal part.
\end{thm}

\medskip
\noindent Analogously to Kurepa's function we consider the
function $A(z)$ given by the integral:
\begin{equation}
\label{A_INT_1}
A(z)
=
\displaystyle\int\limits_{0}^{\infty}{
e^{-t} \displaystyle\frac{t^{z+1}-(-1)^{z}t}{t+1} \: dt},
\end{equation}
which converges for $Re(z) > 0$ \cite{Petojevic:02}, and
it represents one analytical  extension of the alternating sum of
factorials:
\begin{equation}
\label{A_SUM_1}
A(n) = \displaystyle\sum\limits_{i=1}^{n}{(-1)^{n-i}i!} \, .
\end{equation}
For the function $A(z)$ we use term {\em alternating Kurepa's function} and it is one solution of the functional
equation:
\begin{equation}
\label{A_FE_1}
A(z) + A(z-1) = \Gamma(z+1).
\end{equation}
Let us observe that it is possible to make analytical continuation
of alternating Kurepa's function $A(z)$ for $Re(z) \leq 0$.
In that way, the alternating Kurepa's function $A(z)$ is
a meromorphic function with simple poles at $z \!=\! -n$ $(n
\!\geq\! 2)$. Alternating Kurepa's function has the following
residues:
\begin{equation}
\mathop{\mbox{\rm res}}\limits_{z = -n}{A(z)} =
(-1)^n \displaystyle\sum\limits_{k=0}^{n-2}{\displaystyle\frac{1}{k!}}
\quad (n\!\geq\!2).
\end{equation}
Previous results for alternating Kurepa's function are given according to \cite{Petojevic:02}.
The functional equation (\ref{A_FE_1}), besides alternating Kurepa's function $A(z)$, has another
solution by series:
\begin{equation}
\label{Def_A1}
A_{1}(z) = \displaystyle\sum\limits_{n=0}^{\infty}{(-1)^n\Gamma(z+1-n)},
\end{equation}
which converges over the set $\mathbf{C} \backslash \mathbf{Z}$
\cite{Malesevic:06}.

\break

\medskip
\noindent
Extension of domain of functions $A(z)$ and $A_{1}(z)$ in the sense of the principal
part at the point is given by following statements \cite{Malesevic:06}.
\begin{lm}
Let us define $ L_{2} = \displaystyle\sum\limits_{n=0}^{\infty}{
(-1)^{n} \mathop{\mbox{\rm p.p.}}\limits_{z = -(n-1)}{\!\,\Gamma(z)} } $, then$:$
\begin{equation}
\label{Const_L_2}
L_{2}
=
1
+
e \gamma
-
e {\Big (} \displaystyle\sum\limits_{n=1}^{\infty}{\displaystyle\frac{(-1)^{n-1}}{n!n}} {\Big )}
=
1
+
e \mbox{\rm Ei}(-1)
\approx
0.403 \, 652 \, 337 \, ,
\end{equation}
where $\mbox{\rm Ei}$ is function of exponential integral.
\end{lm}
\begin{thm}
For the functions $A(z)$ and $A_{1}(z)$ we have$:$
\begin{equation}
\label{A_PV}
\;\;
\mathop{\mbox{\rm p.p.}}\limits_{z = -n}{\!\!A(z)}
=\! \displaystyle\sum\limits_{i=0}^{n-1}{
(-1)^{n+1-i}\!\!\!\!\mathop{\mbox{\rm p.p.}}\limits_{z=-(i-1)}{\!\!\!\!\Gamma(z)}}
=\! (-1)^{n+1}
{\bigg (}
1\!-\!\displaystyle\sum\limits_{i=1}^{n-1}{\displaystyle\frac{\Gamma^{'}(i)}{\Gamma(i)^2}}
{\bigg )}
\;\; (n \!\in\! \mathbf{N})
\end{equation}
and
\begin{equation}
\mathop{\mbox{\rm p.p.}}\limits_{z = n}{A_{1}(z)}
=
(-1)^n L_{2} + \mathop{\mbox{\rm p.p.}}\limits_{z = n}{A(z)}
\quad (n \!\in\! \mathbf{Z}).
\end{equation}
\end{thm}
The connection between functions $A(z)$ and $A_{1}(z)$ is given by
a formula of the Slavi\'{c}'s type in the following statement
\cite{Malesevic:06}.
\begin{thm}
It is true that$:$
\begin{equation}
\label{A_sum_Slavic}
\;\;\;\;
A(z)
=
{\Big (} e\! \displaystyle\sum\limits_{n=1}^{\infty}{\displaystyle\frac{(-1)^{n-1}}{n! \, n}} - 1 - e \gamma {\Big )}(-1)^z
+
\displaystyle\frac{\pi e}{\sin \pi z}
+
\displaystyle\sum\limits_{n=0}^{\infty}{\!(-1)^{n}\Gamma(z\!+\!1\!-\!n)},
\end{equation}
where the values in the previous formula, in integer points $z$,
are determined in the sense of the principal part.
\end{thm}

\subsection{Principal part of the zeta function and Casimir energy}

Riemann zeta function $\zeta(s)$ is a meromorphic function${}^{\ast)}$\footnote{$\!\!\!\!\!{}^{\ast)}\,$see Example 2.5.
in this paper} and it has only simple pole at $s = 1$ with the principal part:
\begin{equation}
\mathop{\mbox{\rm p.p.}}\limits_{s = 1}{\zeta(s)} = \gamma,
\end{equation}
as cited in \cite{Slavic:70}. We consider the principal part of
the global spectral zeta function $\zeta_{L}(s)$ which is a direct
extension of the Riemann zeta function $\zeta(s)$
\cite{NesterenkoLambiaseScarpetta:05}. Namely, let $L$ be an
elliptic differential operator of the second order acting only on
the variable $x$ and let $\varphi(t,x) = e^{\pm i \omega
t}\varphi_{n}(x)$ be a solution of the following equation:
\begin{equation}
{\bigg (}L + \displaystyle\frac{\partial^2}{c^2 \partial t^2}{\bigg )}\varphi(t,x) = 0,
\end{equation}
where $\omega$ and $c$ are constants. Let scalars $\lambda_{n}$ fulfill
$L \varphi_{n}(x) = \lambda_{n} \varphi_{n}(x)$. Then we define {\em global spectral
zeta function} by \cite{Hawking:77}, \cite{Elizalde:94}, \cite{NesterenkoLambiaseScarpetta:05}:
\begin{equation}
\zeta_{L}(s) = \displaystyle\sum\limits_{n}{\lambda_{n}^{-s}}.
\end{equation}
{\em Casimir energy} of the field $\varphi(t,x)$ is defined by
\cite{BlauVisserWipf:88b}, \cite{NesterenkoLambiaseScarpetta:05}:
\begin{equation}
E_{0}
=
\displaystyle\frac{1}{2}
\lim\limits_{\varepsilon \rightarrow 0}{
\displaystyle\frac{\zeta_{L}(-\frac{1}{2}+\varepsilon) + \zeta_{L}(-\frac{1}{2}-\varepsilon)}{2}}
=
\displaystyle\frac{1}{2}
\mathop{\mbox{\rm p.p.}}\limits_{s = -\frac{1}{2}}{\zeta_{L}(s)}.
\end{equation}
In paper \cite{NesterenkoLambiaseScarpetta:05} some values of
Casimir energy have been given, dependent of the fields which are
considered. All computations in
\cite{NesterenkoLambiaseScarpetta:05}, based on paper
\cite{BlauVisserWipf:88b}, are related to the global spectral zeta
functions with, as a rule, simple poles.


\section{Differential - algebraic properties}

In this section we present a method for proving that certain analytic functions are not
solutions of algebraic  differential equations. The method is based on model-theoretic
properties of differential fields and that $\Gamma(x)$ is a transcendental differential
function.

\subsection{Differential fields}

The theory DF${}_0$ of differential fields of characteristic $0$ is the theory of fields with following axioms
that relate  to the derivative $D$:
\begin{equation}
D(x+y)= Dx + Dy,\quad D(xy)= xDy+yDx.
\end{equation}
Thus, a model of DF${}_0$ is a differential field $\mathbf{K} = (K,+,\cdot, D, 0,1)$ where $(K,+,\cdot,$ $0,1)$
is a field and $D$ is a differential operator satisfying the above axioms. Abraham Robinson proved that DF${}_0$
has a model completion, and then defined DCF${}_0$ to be the model completion of DF${}_0$.
Afterwards Leonore~Blum found simple axioms of DFC${}_0$ not mentioning of differential
polynomials in more than one variable \cite{Sacks:72}. In the following, if not otherwise stated,
$\mathbf{F}, \mathbf{K}, \mathbf{L}, \ldots$ will denote differential fields, $F,L,K,\ldots$
their domains while $\mathbf{F^\ast}, \mathbf{K^\ast}, \mathbf{L^\ast}, \ldots$ will denote their field parts,
i.e. $\mathbf{F^\ast}\mbox{$=$}\mbox{$(F,+,\cdot,0,1)$}$.  It is customary to denote by $\mathbf{L}\{X\}$
the ring of differential polynomials over $\mathbf{L}$ in the variable $X$, see \cite{Marker:96}.
Thus, if $f\in L\{X\}$ then for some natural number $n$, $f=f(X,DX,D^2X,\ldots,D^nX)$ where $f(x,y_1,y_2,\ldots,y_n)$
is the ordinary algebraic polynomial over $\mathbf{L^\ast}$. Then the order of $f$, denoted by $\mbox{\rm ord} \, f$,
is the largest $n$ such that $D^nX$ occurs in $f$. If $f\in L$ we put $\mbox{\rm ord} \, f = -1$ and then we write
$f(a) = f$ for each $a$. For $f \in L\{X\}$ we shall write occasionally $f'$ instead of $Df$, and $f(a)$ instead of
$f(a,Da,D^2a,\ldots,D^na)$ for each $a$ and $n = \mbox{\rm ord} \, f$.

\smallskip
If $b\in K$, then $\mathbf{L}(b)$ will denote the simple differential extension
of $\mathbf{L}$ in $\mathbf{K}$, i.e. $\mathbf{L}(b)$ is the smallest differential
subfield of $\mathbf{K}$ containing both $L$ and $b$. Also, we shall use the
following abbreviations:

\break

\medskip
d.p. is standing for {\em differential polynomial}. Thus, $f$ is a d.p.
over $\mathbf{L}$ in the variable $X$ if and only if $f\in L\{X\}$.

\smallskip
d.a. is standing for {\it differential algebraic}. Hence, if $\mathbf{L} \subseteq \mathbf{K}$
then $b\in K$ is d.a. over $\mathbf{L}$ if and only if there is a non-zero d.p. $f$ such that
$f(b,Db,D^2b,\ldots,D^nb)=0$, otherwise $b$ is transcendental. The field $\mathbf{K}$ is a d.a.
extension of $\mathbf{L}$ if every $b\in K$ is d.a. over $\mathbf{L}$.

\smallskip
d.e. is standing for {\it differential equation},

\smallskip
a.d.e. is standing for {\it algebraic differential equation}. Hence, $f=0$ is a.d.e. if $f\in L\{X\}$.

\medskip
Models of DCF${}_0$ are differentially closed fields. A differential field $\mbox{\bf K}$ is {\em differentially closed}
if, whenever $f, g \in \mbox{\bf K}\{X\}$, $g$ is non-zero and $\mbox{\rm ord} \, f > \mbox{\rm ord} \, g$, there is
$a \in K$ such that $f(a) = 0$ and $g(a) \neq 0$. The theory DCF${}_0$ admits elimination of quantifiers and it is
submodel complete (A.~Robinson): if $\mathbf{F} \subseteq \mathbf{L}, \mathbf{K}$ then $\mathbf{K}_F
\equiv \mathbf{L}_F$, i.e. $\mathbf{K}$ and $\mathbf{L}$ are elementary equivalent over $\mathbf{L}$.
In the following we shall use the next theorem, see \cite{MijajlovicMalesevic:06}:

\begin{thm}
\label{Th_2_0}
Suppose $\mathbf{F} \subseteq \mathbf{K}$ and let $L = \{b \in K \colon \, b \enskip
\mbox{is d.a. over} \enskip \mathbf{F}\}$. Then

\smallskip
\noindent
{\bf a.} $L$ is a differential subfield of $\mathbf{K}$ extending $\mathbf{F}$.

\smallskip
\noindent
{\bf b.} If $\mathbf{K}$ is d.a. closed then $\mathbf{L}$ is d.a. closed.
\end{thm}

Other notations, notions and results concerning differential fields that will
be used corresponds to those in \cite{Sacks:72} or \cite{Marker:96}.

\subsection{Transcendental differential functions}

Suppose $\mathbf{L}\subseteq \mathbf{K}$. Let ${\cal R} = \mathbf{R}(x)$ be the differential
field of real rational functions and ${\cal C} = \mathbf{C}(z)$ the differential field
of complex rational functions. The following H\"older's famous theorem asserts
the differential transcendentality of Gamma function.

\begin{thm}
\label{Th_2_1}
{\bf a.} $\Gamma(x)$ is not d.a. over $\mathbf{R}(x)$.\enskip
{\bf b.} $\Gamma(z)$ is not d.a. over $\mathbf{C}(z)$.\hfill$\stop$
\end{thm}

Now we shall use the transcendentality of $\Gamma(z)$ and properties of differential fields
to prove differential transcendentality over ${\cal C}$ of some analytic functions. Let us
denote by ${\cal M}_{\mathbf{D}}$ the class of complex functions meromorphic on a complex domain
$\mathbf{D}$ (a connected open set in the complex $z$-plane $\mathbf{C}$).
If $\mathbf{D} = \mathbf{C}$ then we shall write ${\cal M}$ instead of ${\cal M}_{\mathbf{D}}$.
Then ${\cal M}_{\mathbf{D}}$ is differential field and ${\cal C} \subseteq {\cal M}$. Further,
let ${\cal L} = \{ f \in {\cal M} \,\mbox{:}\, f \enskip \mbox{d.a. over}\enskip {\cal C}\}$.
By Theorem \ref{Th_2_0} ${\cal L}$ is a differential subfield of ${\cal M}$ extending ${\cal C}$.
The function $\Gamma(z)$ is meromorphic and by H\"older's theorem $\Gamma(z) \not\in {\cal L}$.

\begin{ex}
As we have seen in the first part, Kurepa's function $K(z)$ can be continued meromorphically to whole complex plane.
Therefore, $K(z-1)$ is meromorphic either. Also, as we have seen in the first part, Kurepa's function satisfies
the recurrence relation
\begin{equation}
\label{K_FE_2}
K(z) - K(z-1) = \Gamma(z),
\end{equation}
Now, suppose that $K(z)$ belongs to ${\cal L}$. Then $K(z)$ satisfies an a.d.e.
\begin{equation}
f(z,y,Dy,D^2y,\ldots,D^ny)=0
\end{equation}
where $f(z,y_1,y_2,\ldots,y_n)\in {\cal C}[x,y_1,y_2,\ldots,y_n]$. Then $K(z+1)$ satisfies the a.d.e.
\begin{equation}
f(z+1,y,Dy,D^2y,\ldots,D^ny)=0
\end{equation}
so $K(z+1)$ belongs to ${\cal L}$. As ${\cal L}$ is a field, by (\ref{K_FE_2}) it follows
that $\Gamma(z)$ belongs to ${\cal L}$, what yields a contradiction.
Hence, $K(z)$ is a transcendental differential function.

Using previous method we can conclude that each meromorphic solution of a functional equation
{\rm (\ref{K_FE_2})} is transcendental differential function over the field ${\cal C}$.
For example, another solution of this functional equation is series {\rm (\ref{Def_K1})}.
Therefore, $K_{1}(z)$ is a transcendental differential function too.\hfill$\stop$
\end{ex}

In a similar way, one can prove that alternating functions $A(z)$ and $A_1(z)$
are differentially transcendental too over ${\cal C}$.

\medskip
The following general proposition concerning transcendental differential fu\-ncti\-ons holds.

\begin{thm}
\label{Th_2_3}
Let $a(z)$ be a meromorphic differentially transcendental function over ${\cal C}$
and $f(z,u_0,u_1,\ldots,u_m,y_1,\ldots,y_n)$ be a polynomial over ${\cal C}$.
If $b$ is meromorphic and $f(z,b,E_1b,\ldots,$ $E_mb,Db,\ldots,D^nb)=a(z)$,
where $E_{i} f(z) = \alpha_{i} z + \beta_{i}$, $\alpha_{i}, \beta_{i} \in \mathbf{C}$,
then $b$ is differentially transcendental over ${\cal C}$.
\end{thm}

\noindent
{\bf Proof} \enskip Suppose that $b$ is d.a. over ${\cal C}$, i.e. that $b \in {\cal L}$.
Then $Db,\ldots,D^nb$ belong to ${\cal L}$. Further, there is a.d.e.
$f(z,y,Dy,\ldots,D^ky)=0$ satisfied by $b$, so $E_ib$ satisfies
\begin{equation}
f(\alpha_iz+\beta_i,y,\alpha_1^{-1}Dy,\ldots,\alpha_k^{-1}D^ky)=0,
\end{equation}
i.e. $E_ib\in {\cal L}$, too. Therefore, $g(z,b,E_1b,\ldots,E_mb,Db,\ldots,D^nb)
\in{\cal L}$, so $a(z)$ belongs to ${\cal L}$, a contradiction. \hfill$\stop$

\begin{ex}
The Riemann zeta function defined by
\begin{equation}
\zeta(s)= \sum_{n=1}^\infty n^{-s}= \prod_{p \enskip \mbox{\scriptsize \rm prime}}
{(1-p^{-s})^{-1}}, \quad \mbox{\rm Re}(s)>1,
\end{equation}
is differentially transcendental over ${\cal C}$ $($Hilbert$)$. First we observe that $\zeta(s)$
can be continued meromorphically to whole complex plane with a simple pole at $s=1$
and that $\zeta(s)$ satisfies the well-known functional equation {\rm \cite{Ivic:85}}:
\begin{equation}
\zeta(s)
=
\chi(s)\zeta(1-s),\quad \mbox{\rm where}\quad
\chi(s)
=
\displaystyle\frac{(2\pi)^s}{
2 \Gamma(s) \cos (\mbox{\small $\displaystyle\frac{\pi s}{2}$})}.
\end{equation}
Now, suppose that $\zeta(s)$ is d.a. over ${\cal C}$, i.e. that $\zeta(s)\in {\cal L}$.
Then $\zeta(1-s)$ and $\zeta(s)/\zeta(1-s)$ belong to ${\cal L}$, too, so $\chi(s)$
belongs to ${\cal L}$. The elementary functions $(2\pi)^s$, and $\cos({\pi s\over 2})$
obviously are d.a. over ${\cal C}$ i.e. they belong to ${\cal L}$. As ${\cal L}$ is a field,
it follows that $\Gamma(z)$ belong to ${\cal L}$, too. But this yield a contradiction,
therefore $\zeta(s)$ is differentially transcendental function over ${\cal C}$.
Generally, Dirichlet $L$-series
\begin{equation}
L_{k}(s)= \displaystyle\sum\limits_{n=1}^{\infty}{\kappa_{k}(n) \displaystyle\frac{1}{n^s}}
\quad (k \in Z),
\end{equation}
where $\kappa_{k}(n)$ is Dirichlet character {\rm \cite{IrelandRosen:82}}, is differentially transcendental
function over ${\cal C}$. This follows from a well-known functional equations
\begin{equation}
L_{-k}(s)
=
2^s \pi^{s-1} k^{-s+\frac{1}{2}} \Gamma(1-s) \cos (\mbox{\small $\displaystyle\frac{\pi s}{2}$}) L_{-k}(1-s)
\end{equation}
and
\begin{equation}
L_{+k}(s)
=
2^s \pi^{s-1} k^{-s+\frac{1}{2}} \Gamma(1-s) \sin (\mbox{\small $\displaystyle\frac{\pi s}{2}$}) L_{+k}(1-s) \, .
\end{equation}
Besides Riemann zeta function $\zeta(s)=L_{+1}(s)$, Dirichlet eta function
\begin{equation}
\eta(s)
=
\displaystyle\sum\limits_{n=1}^{\infty}{(-1)^{n+1}\displaystyle\frac{1}{n^s}}
=
(1-2^{1-s})L_{+1}(s)
\end{equation}
and Dirichlet beta function
\begin{equation}
\beta(s)
=
\displaystyle\sum\limits_{n=0}^{\infty}{(-1)^{n}\displaystyle\frac{1}{(2n+1)^s}}
=
L_{-4}(s)
\end{equation}
are transcendental differential functions as examples of Dirichlet series {\rm \cite{Borwein:87}}. \hfill$\stop$
\end{ex}

\begin{ex}
\label{H_1_Function}
The meromorphic function
\begin{equation}
H_{1}(z)= \sum_{n=0}^{\infty}{1\over (n+z)^2}
\end{equation}
is differentially transcendental over ${\cal C}$. Really, $D^2\ln(\Gamma(z))=H_{1}(z)$,
i.e. $\Gamma(z)$ satisfies the a.d.e. $(D^2\Gamma)\Gamma-(D\Gamma)^2-H_{1}\Gamma^2=0$
over $\mbox{$\cal C$}(H_{1})$. Thus, if $H_{1}$ would be d.a. over ${\cal C}$, then by
Theorem {\rm \ref{Th_2_3}}, $\Gamma$ would be too, what yields a contradiction. Hence, $H_{1}(z)$
is differentially transcendental function over ${\cal C}$.
Let us notice the well-known fact that gamma function is the solution of the following
non-algebraic differential equation: $df(z)/dz = \psi(z) f(z)$ for $f(z) = \Gamma(z)$.
\hfill$\stop$
\end{ex}

\begin{re}
We see that in Example {\rm \ref{H_1_Function}} functions $\Gamma(z)$ and $H_{1}(z)$ are differentially
algebraically dependent, i.e. $g(\Gamma,H_{1})=0$, where $g(x,y)= x''x-(x')^2-yx^2$. We do not know
if similar dependencies exist for pairs $(K,\Gamma)$ and $(\zeta,\Gamma)$. It is very likely that
these pairs are in fact differentially transcendental.
\end{re}

\subsection{Operator {\boldmath $\delta$}}

Let $\mathbf{F}_D = (F,+,\cdot,D,0,1)$ be a differential field and $\theta\in F$. We can introduce
a new differential operator $\delta= \theta\cdot D$, i.e. by $\delta(x)= \theta\cdot D(x)$, $x\in F$.
Then $\mathbf{F}_{\delta} = (F,+,\cdot,\delta,0,1)$ becomes a new differential filed.
Let $\overline{\mathbf{F}}_D$ and $\overline{\mathbf{F}}_{\delta}$ denote differential
closures of fields $\mathbf{F}_D$ and $\mathbf{F}_{\delta}$ respectively.

\begin{pr}
Domains of fields  $\overline{\mathbf{F}}_D$ and $\overline{\mathbf{F}}_{\delta}$
are same, i.e.  $\bar F_{D}  = \bar F_{\delta}$.
\end{pr}

\noindent
{\bf Proof} \enskip If $a\in\bar F_\delta$ then $a$ is a solution of an a.d.e.
$\cal E(\delta)$  in respect to the operator $\delta$. We can substitute in this
equation operator $\delta$ with $\theta\cdot D$, and we shall obtain again an a.d.e.
$\cal E'(D)$ but now in respect to $D$. Then $a$ is a solution of this equation, hence
$a\in\bar F_D$. So we proved that $\bar F_{\delta} \subseteq \bar F_{D}$. On the other hand,
$D = \theta^{-1}\delta$, so we may apply a symmetrical argument, hence $a\in\bar F_D$ implies
$a \in \bar F_\delta$, i.e. $\bar F_{D} \subseteq \bar F_{\delta}$. Therefore, we proved
$\bar F_{D}  = \bar F_{\delta}$.
\hfill$\stop$

\medskip
We can ask the natural question if fields $\overline{\mathbf{F}}_D$,
$\overline{\mathbf{F}}_\delta$ are isomorphic. We observe that it is
not necessary $\mathbf{F}_D \cong \mathbf{F}_\delta$. For example, if $\mathbf{F}
= \mathbf{R}(x)$, $D$ is the ordinary differentiation operator and $\delta=xD$, then the equation
$\delta y= y$ has a solution in $\mathbf{F}_\delta$, $y=x$, while the equation $Dy=y$
has no solution in $\mathbf{F}_D$. Hence $\mathbf{F}_{D} \not\cong \mathbf{F}_{\delta}$. Let us remind that
$\tau \colon F_{\delta} \to F_{D}$ is an isomorphism if $\tau$ satisfies:
\begin{equation}
\;\;\;\;
\tau(x+y)=\tau x+\tau y, \;
\tau(xy)=\tau x \,\tau y,\;
\tau(\delta x)= D\tau(x),\;
\tau(0)=0, \tau(1)=1.
\end{equation}

\smallskip
Under some circumstances the isomorphism exists between fields $\overline{\mathbf{F}}_D$,
$\overline{\mathbf{F}}_\delta$, or between certain intermediate fields. For example, the
conditions will be fulfilled if these fields have functional representation and a particular
differential equation has a solution. Let
${\cal L}_D = \{ f \in {\cal M} \,\mbox{:}\, f \enskip \mbox{d.a. over}\enskip {\cal C}
\enskip\mbox{in respect to the ope-}$ $\mbox{rator}\enskip D\}$ and
${\cal L}_{\delta} = \{ f \in {\cal M} \,\mbox{:}\, f \enskip \mbox{d.a. over}\enskip {\cal C}
\enskip \mbox{in respect to the operator}\enskip \delta\}$.

\begin{thm} If $\theta\in {\cal L}_D$ is non-constant and $g$ is a non-constant solution of
 a.d.e. $Dx = \theta \cdot x$ then ${\cal L}_D \cong {\cal L}_{\delta}$.
\end{thm}

\noindent
{\bf Proof} \enskip First, we observe, using the argument as in the proof of the above proposition,
that domains of ${\cal L}_D$ and ${\cal L}_{\delta}$ are same.
Let $\tau\colon {\cal L}_{\delta}\rightarrow {\cal L}_D$
be defined by $\tau (x)= x\circ g$, where $\circ$ is the composition operator.
We see that $g$ is meromorphic and is a.d. over $\cal C$, therefore $g\in {\cal L}_D$.
$\tau$ is well defined since ${\cal L}_D$ is closed under composition. Obviously it satisfies
$\tau(x+y)=\tau x+\tau y,  \enskip \tau(xy)=\tau x \, \tau y$. Further, as $Dg= \theta\circ g$,
\begin{equation}
\;\;\;\;
\tau(\delta x)
\!=\!
(\theta Dx)\circ g
\!=\!
(\theta\circ g)((Dx)\circ g)
\!=\!
Dg((Dx)\circ g)
\!=\!
D(x\circ g)
\!=\!
D(\tau x).
\end{equation}
$\tau\,$ is 1--1 function, since $g$ takes infinitely many values over a bounded region. Therefore,
\mbox{$\tau\colon {\cal L}_D \cong {\cal L}_{\delta}$}. \hfill$\stop$

\medskip
In the case of $\theta=x$, $x$ here denotes a variable (i.e. the polynomial of the degree one),
we can produce an explicit isomorphism  $\tau\colon {\cal L}_D \cong {\cal L}_{\delta}$.
We can define $\tau$ by $\tau: f\to f\circ g$, $f\in {\cal L}_D$, where $g(x)=e^x$.
Observe that this isomorphism corresponds to the transformation $x=e^z$ in the algorithm
of solving of Euler linear differential equations. This observation give us a new, the algebraic
insight into the classical method of solving Euler and similar types (e.g. Legendre linear equation)
of differential equations. We shall give an illustration by example:

\begin{ex}
Solve $x^3y''' + 3x^2y''-2xy'+2y=0$.
\end{ex}

{\em Solution} \enskip This equation is equivalent to
\begin{equation}
(\delta(\delta-1)(\delta-2) +3\delta(\delta-1)-2\delta+2)y = 0,
\end{equation}
i.e. to the equation  $(\delta^3-3\delta+2)y=0$ in $\cal L_\delta$.
The corresponding equation $(D^3-3D+2)y=0$ in $\cal L_D$ has
general solution $c_1h_1+c_2h_2+c_3h_3$, where
$h_1(x)=e^x, h_2(x)=xe^x, h_3(x)=e^{-2x}$. As
$\tau\colon {\cal L}_D \cong {\cal L}_{\delta}$,
and $\tau^{-1}$ is given by
$\tau^{-1}\colon f\to f\circ g^{-1}$ (here $g^{-1}(x)= \ln x$),
it follows that
\begin{equation}
\tau^{-1}(c_1h_1+c_2h_2+c_3h_3)
=
c_1h_1\circ g^{-1} + c_2h_2\circ g^{-1} + c_3h_3\circ g^{-1}
\end{equation}
is the general solution of $(\delta^3-3\delta+2)y=0$ in $\cal L_\delta$,
and so the solution of the starting equation is $y=c_1x+c_2x\ln x+c_3x^{-2}$.

\medskip
Hence, one should expect that standard methods of solving differential equations which are done by
"properly chosen transformations of the independent variable" correspond in fact to constructions
of an isomorphism between $\cal L_D$ and $\cal L_\delta$, or some other intermediate fields, for
properly chosen differential operators $\delta$.

\end{document}